\documentclass[12pt]{amsart}
\usepackage{amssymb,latexsym, amsmath, amsxtra}
\usepackage{graphicx}
\begin{document}

\title{Moments of zeta and correlations of divisor-sums: I}
 
\author{Brian Conrey}
\address{American Institute of Mathematics, 360 Portage Ave, Palo Alto, CA 94306, USA and School of Mathematics, University of Bristol, Bristol BS8 1TW, UK}
\email{conrey@aimath.org}
\author{Jonathan P. Keating}
\address{School of Mathematics, University of Bristol, Bristol BS8 1TW, UK}
\email{j.p.keating@bristol.ac.uk}

\thanks{We gratefully acknowledge support under EPSRC Programme Grant EP/K034383/1
LMF: L-Functions and Modular Forms.  Research of the first author was also supported by the American Institute of Mathematics and by a grant from the National Science Foundation. JPK is grateful for the following additional support: a grant from the Leverhulme Trust, a Royal Society Wolfson Research 
Merit Award, a Royal Society Leverhulme Senior Research Fellowship, and a grant from the Air Force Office of Scientific Research, Air Force Material Command, 
USAF (number FA8655-10-1-3088). He is also pleased to thank the American Institute of Mathematics for hospitality during a visit where this work started.}

\date{\today}

\begin{abstract}
We examine the calculation of the second  and fourth moments and shifted moments of the Riemann zeta-function on the critical line using long Dirichlet 
polynomials and divisor correlations.  Previously this approach has proved unsuccessful in computing moments beyond the eighth, even heuristically.  A 
careful analysis of the second and fourth moments illustrates the nature of the problem and enables us to identify the terms that are missed in the standard application of these methods.
\end{abstract}

\maketitle

\section{Introduction}
We  revisit the long-standing problem of determining the moments of the Riemann zeta-function on the critical line from a somewhat new perspective. On one hand, we have the detailed conjectures
of [KS, CFKRS] which predict very precisely the outcome of any moment calculation. On the other hand, we have the conventional 
approaches of analytic number theorists via Dirichlet polynomial approximations to $\zeta(s)$ or $\zeta(s)^k$ which proceed
by refinements of the approximate Parseval formula for Dirichlet polynomials (i.e. the Montgomery-Vaughan formula). 
We would like to bring these approaches closer together. In particular, the Dirichlet polynomial approach completely fails 
when considering the 10th or higher moment of $\zeta(s)$ on the critical line: it predicts negative values when the moments are clearly non-negative.  We want to understand 
and rectify, at least from a heuristic perspective, this failing.

Our aim here is to illustrate the nature and cause of this failing by a careful analysis of the second and fourth moments.  We take the standard approach of approximating the 
zeta function by long Dirichlet polynomials and computing moments using divisor correlations, but depart in using longer polynomials than is normal.  

The approximate functional equation expresses the zeta 
function at a height $T$ on the critical line by Dirichlet polynomials of length of the order of $T^{1/2}$.  This allows all principal terms (called, in [CFKRS], the moment polynomials) in the second [I, A] and fourth [HB] moments 
to be determined.  This approach extends heuristically to the sixth moment, and to calculating the leading-order in the asymptotics of the eighth moment; see, for example, [CG].  For moments beyond the eighth, it leads to answers that are clearly incorrect, 
in that they are negative. In order to understand this problem better, we here examine the calculations for the second and fourth moments when we use longer Dirichlet polynomials; specifically, we take the length to be substantially larger than $T$ (and so do not utlize the approximate functional equation).  
We find that in this case the success of the standard methods when applied to the second moment rests on a particular identity, and that they fail when applied to the fourth moment by missing certain key terms.  We believe that this sheds new light on the moment problem.
 
 \section{Second moment}

We begin with the simplest case, that of the second moment of $\zeta(s)$; even in this case our
perspective reveals some new information. 
Our starting point is the theorem of Ingham [I] that
\begin{eqnarray*}
\int_0^T \zeta(s+\alpha)\zeta(1-s+\beta) ~dt =\int_0^T \left( \zeta(1+\alpha+\beta)+\left(\frac t {2\pi}\right)^{-\alpha-\beta}
\zeta(1-\alpha-\beta)\right) ~dt +O(T^{1/2+\epsilon})
\end{eqnarray*}
uniformly for small $\alpha$ and $\beta$, where $s=1/2+it$. Sandro Bettin [B] has proven that, remarkably, the same formula holds 
for a much larger range of $\alpha$ and $\beta$ namely $|\Re \alpha| |\Re \beta| <1/4$
and $|\Im \alpha|, |\Im \beta| \ll T$. 

We  compare this result with what one gets when using the work of Goldston-Gonek  [GG] 
on the mean-square of long Dirichlet polynomials. 
We refer to 
Corollary 1 of [GG] together with example 1, which is about the case $\alpha=\beta=0$ of Ingham's theorem. 
In particular we  approximate 
$\zeta(s+\alpha)$ by a Dirichlet polynomial $\sum_{n\le X} n^{-s-\alpha}$ with $X$ substantially larger than $T$ , say $X=T^\theta$ with 
$1<\theta < 2$; the error in such an approximation is $\ll X^{1/2}/T$ when $t\approx T$
so that $\zeta(s)$ is well-approximated point-wise by this long Dirichlet polynomial.
We introduce a smooth weight function $\psi(u)$ which is real and compactly supported, say on $[1,2]$, 
and then evaluate the mean square of $\zeta$  by evaluating  (as [GG] do) 
\begin{eqnarray*}\int_0^\infty \psi(t/T)|\zeta(1/2+it)|^2 ~dt=
\int_0^\infty \psi(t/T)\left|\sum_{n\le X}\frac{1}{n^{1/2+it}}\right|^2 ~dt+O(X^{1/2}\log T).
\end{eqnarray*}
Using coefficient correlations, Goldston and Gonek prove that  this   is
\begin{eqnarray*}
= \hat{\psi}(0)T \sum_{n\le X}\frac 1 n +2 \Re T \int_y^\infty \big [ \frac{v}{y}\big] \hat{\psi}(v) \frac{dv}{v} +O(T^\epsilon(X^{1/2} +X/T))
\end{eqnarray*}
where 
$$y= \frac{T}{2\pi X},$$
 $[v]$ is the floor of $v$, i.e.~the greatest integer less than 
or equal to $v$, and the Fourier transform is defined by
$$\hat{\psi}(v)=\int_{\mathbb R} \psi(u) e(uv)~du$$
where $e(x)=\exp(2\pi i x)$.

On the other hand, 
the mean square of $\zeta$ is [A] 
\begin{eqnarray*}
T\int_0^\infty \psi(t) \left(\log\frac{tT}{2\pi } +2\gamma\right) ~dt +O(X^{1/2}).
\end{eqnarray*}
Using the fact that
$$\sum_{n\le X} \frac{1}{n} =\log X +\gamma +O(1/X)$$
it follows that  
\begin{eqnarray} \label{eqn:floor}
2\Re \int_y^\infty \big [ \frac{v}{y}\big] \hat{\psi}(v) \frac{dv}{v}
\sim \int_0^\infty \left(\log uy+\gamma\right) \psi(u) ~du
\end{eqnarray}
as $y\to 0^+$.  In the present formulation, this is purely a statement about a relationship between a function and its 
Fourier transform for real, compactly supported functions $\psi$.  

We note that establishing the compatibility between the result of [GG] and that in [I, A] in the  approach taken here relies on proving equation (\ref{eqn:floor}).  This 
we do below, as a consequence of a more general identity.

\subsection{General $\alpha $ and $\beta$} If we do this with small shifts $\alpha$ and $\beta$ we have, again with $s=1/2+it$, 
\begin{eqnarray*}
\int_{\mathbb R} \psi\left( \frac t T \right) \zeta(s+\alpha)\zeta(1-s+\beta) ~dt
&=&T\int_{\mathbb R} \psi(t) \left( \zeta(1+\alpha+\beta)+\left(\frac {tT}{2\pi}\right) ^{-\alpha-\beta}
\zeta(1-\alpha-\beta)\right) ~dt\\
&&\qquad  +O(T^{1/2+\epsilon})
\end{eqnarray*}
by Ingham. By [GG] this is also equal to 
\begin{eqnarray*}
\hat \psi(0)T \sum_{n\le X} \frac{1}{n^{1+\alpha+\beta}}
+2T  \int_y^\infty \sum_{h\le \frac v y}
\left(\frac{ hT }{2\pi v}\right)^{-\alpha-\beta}\Re \hat{\psi}(v)
\frac{dv}{v}+O(T^\epsilon(X/T +X^{1/2}))
\end{eqnarray*}
where $y=\frac{T}{2\pi X}$.
Since 
$$\sum_{n\le U}\frac{1}{n^{1+\alpha+\beta}}=\zeta(1+\alpha+\beta) -\frac{U^{-\alpha-\beta}}{\alpha+\beta}+O(1/U)$$
for $\Re (\alpha+\beta)>0$ and $U> T/\pi$, it follows that 
\begin{eqnarray*}
\int_{\mathbb R} \psi(t) \left(\frac {y^{\alpha+\beta}}{\alpha+\beta}+
t ^{-\alpha-\beta}\zeta(1-\alpha-\beta)\right) ~dt
\sim 2\Re \int_y^\infty \sum_{h\le \frac v y } 
h^{-\alpha-\beta}\hat{\psi}(v)
\frac{dv}{v^{1-\alpha-\beta}}
\end{eqnarray*}
as $y\to 0^+$. 
In fact,   this is an identity. 
Note that (\ref{eqn:floor}) is the limiting case when $\alpha, \beta \to 0$ of this identity.

When we pass to higher moments the corresponding identities will be consequences of 
conjectures rather than facts. Thus, it is desirable to give a direct proof of this now, one that 
can be imitated later. 

\begin{proof}
Write $\delta= \alpha+\beta$. We can extend the integration in $v$ down to $v=0$ since the integrand is 0 for $0<v<y$. 
 We will show that
\begin{eqnarray} \label{eqn:id}
\int_{\mathbb R} \psi(t) \left(\frac {y^\delta}{\delta}+
t ^{-\delta}\zeta(1-\delta)\right) ~dt
=2\Re \int_0^\infty \sum_{h\le \frac v y } 
h^{-\delta}\hat{\psi}(v)
\frac{dv}{v^{1-\delta}} .
\end{eqnarray} 
To begin with use
$$\sum_{h\le \frac v y}h^{-\delta} =\frac{1}{2\pi i} \int_{(2)} \zeta(s+\delta) \frac{ (\frac vy  )^s}{s}~ds.$$
Then write
$$2\Re \hat{\psi}(v)=\int_{\mathbb R} \psi(u) (e(-uv)+e(uv)) ~du.$$
Because of the support of $\psi$ the $u$-integral only runs over a finite segment of 
the positive reals. We now want to 
 interchange the $u$-integration and the $v$- integration, but we have to be careful here. Basically
we first split the $u$-integrand into two pieces, one with $e(uv)$ and one with $e(-uv)$ and we deform the 
paths of the two ensuing $v$-integrals by rotating the straight line paths along the positive reals onto 
the positive imaginary axis (in the case of the $e(uv)$ integral) and  the negative imaginary axis (in the case of the 
$e(-uv)$ integral). In this way we get that the right hand side of (\ref{eqn:id}) is 
\begin{eqnarray*}
\int_0^\infty \psi (u) \frac{1}{2\pi i}\int_{(2)} \zeta(s+\delta) \frac{y^{-s}}{s}
\left(\int_0^\infty (iv)^{s+\delta} e^{-2\pi v u} \frac{dv}{v} +\int_0^\infty (-iv)^{s+\delta} e^{-2\pi v u} \frac{dv}{v} 
\right) ~ds ~du.
\end{eqnarray*}
The sum of the two $v$-integrals is 
$$2\cos \frac{\pi (s+\delta)}{2} (2\pi)^{-s-\delta} \Gamma(s+\delta) u^{-s-\delta}$$ which is 
$u^{-s-\delta}\chi(1-s-\delta)$ where $\chi$ is the factor from the functional equation $\zeta(1-s)=\chi(1-s)\zeta(s)$. 
Thus, our expression simplifies to 
\begin{eqnarray} \label{eqn:mellin} 
\int_0^\infty \psi (u) \frac{1}{2\pi i}\int_{(2)} \zeta(1-s-\delta) \frac{y^{-s}u^{-s-\delta}}{s}
  ~ds ~du.
\end{eqnarray}
Now we move the $s$-path of integration to the left to $\Re s=-\infty$. When we account for the residues from the poles at 
  $s=0$ and $s=-\delta$ we exactly obtain the left-hand side of (\ref{eqn:id}). 
\end{proof}

Note that we have shown that if $\Re s>0$ then  
\begin{eqnarray} \label{eqn:chi}
2\int_0^\infty v^s ~ \Re \hat \psi(v) \frac{dv}{v}
= \chi(1-s)\int_0^\infty \psi(t) t^{-s} ~dt
\end{eqnarray}
for functions $\psi $ compactly supported on $(0,\infty)$.
 
We conclude that one can calculate the second moment of the zeta function
with a power savings on the error term  using long Dirichlet polynomials (i.e.~longer than in the 
approximate functional equation), but this relies on establishing certain 
identities, such as that proved above.  We shall see that this becomes more challenging for the higher moments.

\section{Fourth moment}
The fourth moment of the Riemann zeta-function has been much studied
with the explicit formula of Motohashi [M] one of the crowning achievements.

The main goal of this paper is to  revisit this topic but through the lens described in the last section.
  We want 
to understand how one can approach the fourth moment via the mean square 
of a long Dirichlet polynomial approximation to $\zeta(s)^2$. This point of view
has also been explored to a certain extent in the work [H-B] of Heath-Brown,
which built upon that of Ingam and Atkinson [I, A]; see also [CG].
But each of the above used an approximate functional equation 
to reduce the necessary length of the approximating 
polynomial. Here we want to see what happens purely with Dirichlet
polynomials. This perspective will turn out to be 
useful in a future discussion of yet higher moments and divisor correlations. 

\section{The 4th moment of zeta}
To begin with we recall the fourth moment of zeta. Let $s=1/2+it$
and suppose that $\alpha, \beta, \gamma, \delta \ll (\log T)^{-1}$.
Then 
\begin{eqnarray} \label{eqn:fourth}&&
\int_0^\infty \psi\left(\frac tT\right) \zeta(s+\alpha)\zeta(s+\beta)\zeta(1-s+\gamma)\zeta(1-s+\delta)~dt\\
&&\qquad \qquad \qquad \nonumber
=T\int_0^\infty \psi(t) \mathcal Z_{tT}(\alpha,\beta,\gamma,\delta) ~dt +O(T^{2/3+\epsilon}),
\end{eqnarray}
where
\begin{eqnarray} \label{eqn:swaps}
\mathcal Z_t(\alpha,\beta,\gamma,\delta)&=&  \nonumber
Z(\alpha,\beta,\gamma,\delta) +\left(\frac{t}{2\pi}\right)^{-\alpha-\gamma}Z(-\gamma,\beta,-\alpha,\delta)
+\left(\frac{t}{2\pi}\right)^{-\alpha-\delta}Z(-\delta,\beta,\gamma,-\alpha)\\
&&\qquad +\left(\frac{t}{2\pi}\right)^{-\beta-\gamma}Z(\alpha,-\gamma, -\beta,\delta)
+\left(\frac{t}{2\pi}\right)^{-\beta-\delta}Z(\alpha,-\delta, \gamma, -\beta)\\
&&\qquad \qquad \nonumber
+\left(\frac{t}{2\pi}\right)^{-\alpha-\beta -\gamma-\delta }Z(-\gamma,-\delta, -\alpha, -\beta )\end{eqnarray}
where 
\begin{eqnarray*}
Z(\alpha,\beta,\gamma, \delta)=\frac{\zeta(1+\alpha+\gamma)\zeta(1+\alpha+\delta)\zeta(1+\beta+\gamma)\zeta(1+\delta)}
{\zeta(2+\alpha+\beta+\gamma+\delta)}
\end{eqnarray*}
 (Perhaps ``recall'' is not the correct word here since the above is certainly a theorem but its proof is 
not written down in full details anywhere.)

Now
$$\zeta(s+\alpha)\zeta(s+\beta) =\sum_{n=1}^\infty \frac{\tau_{\alpha,\beta}(n)}{n^{2s}}$$
where
$$\tau_{\alpha,\beta}(n):=\sum_{hk=n} h^{-\alpha}k^{-\beta}.$$
We let 
$$D_{\alpha,\beta}(s)=\sum_{n=1}^\infty \frac{\tau_{\alpha,\beta}(n)}{n^s}$$
and 
$$D_{\alpha,\beta}(s;X)=\sum_{n\le X} \frac{\tau_{\alpha,\beta}(n)}{n^s}$$
be its approximating Dirichlet polynomial of length $X$.
Let
$$I(T;X)=I_{\alpha,\beta,\gamma,\delta}(T;X)=\int_0^\infty \psi\left(\frac t T\right) D_{\alpha,\beta}(s,X)D_{\gamma,\delta}(1-s,X) ~dt.$$
Our basic question is: How does $I(T;X)$ approach 
$\int_0^\infty \psi\left(\frac t T\right) \zeta(s+\alpha)\zeta(s+\beta)\zeta(1-s+\gamma)\zeta(1-s+\delta)~dt$ as $X\to \infty$?

\section{Descending through the recipe}
If we cheat, we can work backwards in a way and deduce the behavior of  $I(T;X)$ from that of 
$\int_0^\infty \psi\left(\frac t T\right) \zeta(s+\alpha)\zeta(s+\beta)\zeta(1-s+\gamma)\zeta(1-s+\delta)~dt$.
The starting point is that by Perron's formula
$$D_{\alpha,\beta}(s;X)=\frac{1}{2\pi i}\int_{(2)} D_{\alpha,\beta}(s+z) X^z \frac{dz}{z}.$$
So, then
$$I(T;X)= \frac{1}{(2\pi i)^2 }\iint_{z,w} X^{w+z} \int_0^\infty \psi\left(\frac t T\right) D_{\alpha,\beta}(s+z)D_{\gamma,\delta}(1-s+w)
~dt \frac{dz}{z}\frac {dw}{w}.$$
Now, if we move the paths of integration in $z$ and $w$ so that their real parts are positive but small, then we can
use the formula (\ref{eqn:fourth}) to evaluate the integral over $t$; however the imaginary parts of $z$ and $w$ 
can be large. We believe (\ref{eqn:fourth}) holds uniformly when the real parts of the shifts 
are $\ll 1/\log T$ but the imaginary 
parts of the shifts can be large, as large as $T^{1-\epsilon}$ and even larger. This has not been proven
to hold for (\ref{eqn:fourth}), but has been proven, as mentioned above, by Sandro Bettin [B] in the case of the second moment. 

Now we use the fact that 
$$D_{\alpha,\beta}(s+z)=D_{\alpha+z,\beta+z}(s)$$ 
and, assuming the uniform version of (\ref{eqn:fourth}), replace 
 $\int_0^\infty \psi\left(\frac t T\right) D_{\alpha,\beta}(s+z)D_{\gamma,\delta}(1-s+w)
~dt$ in the above by 
$$T\int_0^\infty \psi(t) \mathcal Z_{tT}(\alpha+z, \beta+z, \gamma+w, \delta+w) ~dt.$$
Next, note that 
$$\mathcal Z_t(\alpha+z, \beta+z, \gamma+w, \delta+w)=\mathcal Z_t(\alpha+z+w, \beta+z+w, \gamma, \delta)$$
follows just from the symmetries of $\mathcal Z_t$. 
We now change variables with $s=z+w$ and have 
$$I(T;X)= \frac{1}{(2\pi i)^2 }\int_{\Re s =4}\int_{\Re w=2} \frac{X^{s}}{w (s-w)}T\int_0^\infty \psi(t) \mathcal Z_{tT}(\alpha+s,\beta+s;\gamma,\delta)~dt
~dw~ds +O(T^{2/3+\epsilon}) 
 .$$
We move the path of integration in the $w$-variable off to the left, to $-\infty$ and collect the residue from
the pole at $w=0$. In this way we have 
$$I(T;X)= \frac{1}{2\pi i }\int_{\Re s =4}  \frac{X^{s}}{s} T\int_0^\infty \psi(t) \mathcal Z_{tT}(\alpha+s,\beta+s;\gamma,\delta)~dt
~ds +O(T^{2/3+\epsilon}) 
 .$$

\section{How many swaps to use?}
In the formula (\ref{eqn:swaps}) we sort the terms according to the number $U$ of ``swaps''; this is the 
number of pairs of variables that are exchanged in the arguments of $\mathcal Z_t$. Specifically, there are 6 terms on the right side
of (\ref{eqn:swaps}), the first has 0-swaps, the next four each have one swap, and the sixth and last term has two swaps.
This notion is relevant because in the expressions for 
$$ \mathcal Z_t(\alpha+s,\beta+s;\gamma,\delta)$$ the terms with $U$ swaps have  
factors of the form $\left(\frac {tT} {2\pi}\right)^{-Us}.$
When paired with the factor $X^s$ we see that to evaluate the relevant integral over $s$ we should
move the path to the right (and get 0 for an answer) whenever $X< \left(\frac {tT} {2\pi}\right)^{U};$
otherwise we should move the path to the left and evaluate the integral as the sum 
of the residues over the poles.  

\section{Explicit formula: 0-swaps}
We refer to the terms with 0-swaps as the ``diagonal'' terms. These are just
$$
\hat{\psi}(0)T\frac{1}{2\pi i }\int_{\Re s =4}  \frac{X^{s}}{s}  \mathcal Z(\alpha+s,\beta+s;\gamma,\delta)
~ds  
 $$
which lead to $T\hat{\psi}(0)$ times 
\begin{eqnarray*}&&
\frac{1}{2\pi i}\int_{(2)}\frac{\zeta(s+1+\alpha+\gamma)\zeta(s+1+\alpha+\delta)\zeta(s+1+\beta+\gamma)\zeta(s+1+\beta+\delta)}
{\zeta(2s+2+\alpha+\beta+\gamma+\delta)}\frac{X^s}{s}~ds\\
&&\qquad = \frac{\zeta(1+\alpha+\gamma)\zeta(1+\alpha+\delta)\zeta(1+\beta+\gamma)\zeta(1+\beta+\delta)}
{\zeta(2+\alpha+\beta+\gamma+\delta)}\\
&&\qquad \qquad
 - \frac{X^{-\alpha-\gamma}}{\alpha+\gamma}\frac{\zeta(1+\delta-\gamma)\zeta(1+\beta-\alpha)
\zeta(1+\beta-\alpha+\delta-\gamma)}{\zeta(2+\beta-\alpha+\delta-\gamma)}\\
&&\qquad  \qquad \qquad- \frac{X^{-\alpha-\delta}}{\alpha+\delta}\frac{\zeta(1+\gamma-\delta)\zeta(1+\beta-\alpha)
\zeta(1+\beta-\alpha+\gamma-\delta)}{\zeta(2+\beta-\alpha+\gamma-\delta)}\\
&&\qquad  \qquad \qquad \qquad - \frac{X^{-\beta-\gamma}}{\beta+\gamma}\frac{\zeta(1+\beta-\alpha)\zeta(1+\delta-\gamma)
\zeta(1+\alpha-\beta+\delta-\gamma)}{\zeta(2+\alpha-\beta+\delta-\gamma)} \\
&& \qquad \qquad \qquad \qquad \qquad -  \frac{X^{-\beta-\delta}}{\beta+\delta}\frac{\zeta(1+\alpha-\beta)\zeta(1+\gamma-\delta)
\zeta(1+\alpha-\beta+\gamma-\delta)}{\zeta(2+\alpha-\beta+\gamma-\delta)}
\end{eqnarray*}
plus an error term which is  $O(T^{2/3+\epsilon})$.

\section{Explicit formula: one swap}
There are four terms with one-swap. One of these is 
\begin{eqnarray} \label{eqn:1swap}
T\int_0^\infty \psi(t)  \frac{1}{2\pi i }\int_{\Re s =4}  \frac{X^{s}}{s}  \left(\frac{tT}{2\pi}\right)^{-\alpha-\gamma-s}Z(-\gamma,\beta+s,-\alpha-s,\delta)
~ds ~dt.
\end{eqnarray}
If $X<\frac{T}{2\pi}$ this is 0. If $X>\frac{T}{\pi}$ then we evaluate it as a sum over the
residues at all of the poles (near $s=0$) of the integrand.
Note that
$$Z(-\gamma,\beta+s,-\alpha-s,\delta)=\frac{\zeta(1-\alpha-\gamma-s)\zeta(1-\gamma+\delta)\zeta(1-\alpha+\beta)\zeta(1+\beta+\delta+s)}
{\zeta(2-\alpha+\beta-\gamma+\delta)}
$$
which has poles at $s=-\alpha-\gamma$ and $s=-\beta-\delta$.
 The sum of the residues at these poles and at $s=0$   give
\begin{eqnarray} \label{eqn:oneswap}&&
T\int_0^\infty \psi(t) 
 \bigg(
\left(\frac {tT}{2\pi}\right)^{-\alpha-\gamma} \frac{\zeta(1-\alpha-\gamma) \zeta(1-\gamma+\delta)\zeta(1-\alpha+\beta)\zeta(1+\beta+\delta)}
{\zeta(2-\alpha+\beta-\gamma+\delta)}\\
&&\qquad \nonumber +
\frac{X^{-\alpha-\gamma}}{\alpha+\gamma}
  \frac{ \zeta(1-\gamma+\delta)\zeta(1-\alpha+\beta)\zeta(1-\alpha +\beta-\gamma+\delta)}
{\zeta(2-\alpha+\beta-\gamma+\delta)}\\
&& \qquad \qquad \nonumber  -
\frac{X^{-\beta-\delta}}{\beta+\delta}  \left(\frac {tT}{2\pi}\right)^{-\alpha+\beta-\gamma+\delta} \frac{\zeta(1-\alpha+\beta-\gamma+\delta)
\zeta(1-\gamma+\delta)\zeta(1-\alpha+\beta) }
{\zeta(2-\alpha+\beta-\gamma+\delta)}\bigg) ~dt.
\end{eqnarray}
Note that the term in the integrand which is independent of $t$ cancels with a term from the 0-swap diagonal piece above. 

\section{Explicit formula: two swaps}
For the term with two-swaps we have
\begin{eqnarray*}
T\int_0^\infty \psi(t)  
\frac{1}{2\pi i}\int_{(2)} \left(\frac{tT}{2\pi}\right)^{-\alpha-\beta -\gamma-\delta -2s}Z(-\gamma,-\delta, -\alpha-s, -\beta-s )
\frac{X^s}{s} ~ds ~dt .\end{eqnarray*}
If $T^2>\pi^2 X$ then the answer will be 0 as can be seen by moving the path of integration to the right.
If  $T^2< 4\pi^2 X$ then we use the fact that
\begin{eqnarray*}
Z(-\gamma,-\delta, -\alpha-s, -\beta-s )=\frac{\zeta(1-\alpha-\gamma-s)\zeta(1-\beta-\gamma-s)\zeta(1-\alpha-\delta-s)\zeta(1-\beta-\delta-s)}
{\zeta(2-\alpha-\beta-\gamma-\delta-2s)}
\end{eqnarray*}
has poles at $s=\alpha+\gamma, \alpha+\delta, \beta+\gamma$ and $\beta+\delta$;
 so together with the residue from the pole at $s=0$ we will get a total of five terms 
for this integral:
\begin{eqnarray*}&&
\left(\frac{tT}{2\pi}\right)^{-\alpha-\beta -\gamma-\delta}
\frac{\zeta(1-\alpha-\gamma)\zeta(1-\beta-\gamma)\zeta(1-\alpha-\delta)\zeta(1-\beta-\delta)}
{\zeta(2-\alpha-\beta-\gamma-\delta)}\\
&&\qquad + \frac{X^{-\alpha-\gamma}}{\alpha+\gamma}\left(\frac{tT}{2\pi}\right)^{\alpha-\beta +\gamma-\delta}
\frac{ \zeta(1+\alpha -\beta)\zeta(1+\gamma-\delta)\zeta(1+\alpha-\beta+\gamma-\delta)}
{\zeta(2+\alpha-\beta+\gamma-\delta)}\\
&&\qquad + \frac{X^{-\alpha-\delta}}{\alpha+\delta}\left(\frac{tT}{2\pi}\right)^{\alpha-\beta -\gamma+\delta}
\frac{ \zeta(1-\gamma+\delta)\zeta(1+\alpha-\beta-\gamma+\delta)\zeta(1+\alpha-\beta )}
{\zeta(2+\alpha-\beta-\gamma+\delta)}\\
&&\qquad + \frac{X^{-\beta-\gamma}}{\alpha+\delta}\left(\frac{tT}{2\pi}\right)^{-\alpha+\beta +\gamma-\delta}
\frac{ \zeta(1-\alpha+\beta)\zeta(1-\alpha+\beta+\gamma-\delta)\zeta(1+\gamma-\delta)}
{\zeta(2-\alpha+\beta+\gamma-\delta)}\\
&&\qquad + \frac{X^{-\beta-\delta}}{\beta+\delta}\left(\frac{tT}{2\pi}\right)^{-\alpha+\beta -\gamma+\delta}
\frac{ \zeta(1-\alpha+\beta-\gamma+\delta)\zeta(1- \gamma+\delta)\zeta(1-\alpha+\beta)}
{\zeta(2-\alpha+\beta-\gamma+\delta)}
\end{eqnarray*}
Notice that the last term here is the negative of the last term in the explicit formula for one-swap. 
In the case that $X>\left(\frac T{2\pi}\right)^2$ we add the 0-swap, one swap, and two swap terms all together 
and all of the dependency on $X$ disappears; we recover exactly the formula (\ref{eqn:fourth}). This is to be expected
since when $X>\left(\frac T{2\pi}\right)^2$ it is known that the  Dirichlet polynomial $D_{\alpha,\beta}(s;X) $ is a good point-wise
approximation to $\zeta(s+\alpha)\zeta(s+\beta)$, so the mean squares are the same.

\section{Ascending through convolution sums: the diagonal}
Next we turn to what can be said about the mean-square $I(T;X)$ but now from the point 
of view of coefficient correlations. We call this the ``ascending'' perspective.
 We consider
$$I ^\psi(T,X)=
\int_0^\infty  D_{\alpha,\beta}(s,X)D_{\gamma,\delta}(1-s,X) \psi(t/T) ~dt
=T\sum_{m,n\le X} \frac{\tau_{\alpha,\beta}(m)\tau_{\gamma,\delta}(n)}{\sqrt{mn}} \hat \psi(T\log(m/n).$$
The simplest terms are the diagonal terms:
$$\hat{\psi}(0) T \sum_{n\le X} \frac{\tau_{\alpha,\beta}(n)\tau_{\gamma,\delta}(n)}{n}.$$

The sum here may be evaluated by Perron's formula as
$$\sum_{n\le X}\frac{\tau_{\alpha,\beta}(n)\tau_{\gamma,\delta}(n)}{n}=
\frac{1}{2\pi i}\int_{(2)}\sum_{n=1}^\infty  \frac{\tau_{\alpha,\beta}(n)\tau_{\gamma,\delta}(n)}{n^{1+s}}\frac{X^s}{s} ~ds.
$$
But 
$$\sum_{n=1}^\infty  \frac{\tau_{\alpha,\beta}(n)\tau_{\gamma,\delta}(n)}{n^{1+s}}
=\mathcal Z(\alpha+s,\beta+s;\gamma,\delta)
$$
and so the diagonal terms are identical to the 0-swap terms already considered above.
 
\section{Ascending: Type I off-diagonals with shifts}
Now we are left with the terms
$$
=T\sum_{m,n\le X\atop m\ne n} \frac{\tau_{\alpha,\beta}(m)\tau_{\gamma,\delta}(n)}{\sqrt{mn}} \hat \psi(T\log(m/n)).$$
We let $m=n+h$ and obtain (just as in [GG])
\begin{eqnarray*}
2T \Re \int_y^\infty  
\sum_{h\le \frac v y } 
M'\left(\frac {hT}{2\pi v},h\right)  \hat{\psi}(v)
\frac{dv}{v} 
\end{eqnarray*}
where the coefficient correlation is
$$\sum_{n\le u}\tau_{\alpha,\beta}(n)\tau_{\gamma,\delta}(n+h)=M(u,h)+E(u,h).$$

 We 
use the $\delta$-method of [DFI] to get an expression for $M'(u,h)$,   the average value of  
$\sum_{n\le u}  \tau_{\alpha,\beta}(n) \tau_{\gamma,\delta}(n+h)$. Namely,
$$ \sum_{q=1}^\infty  c_q(h) 
P_{\alpha,\beta}(u,q)P_{\gamma,\delta}(u+h,q)$$
where $c_q(h)=\sum_{d\mid h\atop d\mid q} d\mu(q/d)$ is the Ramanujan sum and 
where  $P_{\alpha,\beta}(u,d)$ is 
the average value of $\sum_{n\le u} \tau_{\alpha,\beta}(n)e(n/d)$. The latter is just the sum of the residues near 1 of 
$\mathcal D_{\alpha,\beta}(s,1/d)u^{s-1}$ where 
$$\mathcal D_{\alpha,\beta}(s,1/d) =\sum_{n=1}^\infty \tau_{\alpha,\beta}(n)e(n/d) n^{-s}.$$ 
In [C] Lemma 4, for example,  it is shown that $\mathcal D_{\alpha,\beta}(s,1/q)-q^{1-2s-\alpha-\beta}D_{\alpha,\beta}(s)$ is entire.
Thus,
\begin{eqnarray*}
P_{\alpha,\beta}(u,q) = \frac 1 q \zeta(1-\alpha+\beta) u^{-\alpha}q^{\alpha-\beta}
+\frac 1 q \zeta(1+\alpha-\beta) u^{-\beta}q^{\beta-\alpha}  
\end{eqnarray*}
and
\begin{eqnarray*}&&
M'(u,h) = \sum_{q=1}^\infty \frac{c_q(h)}{q^2}\left(\zeta(1-\alpha+\beta) u^{-\alpha}q^{\alpha-\beta}
+ \zeta(1+\alpha-\beta) u^{-\beta}q^{\beta-\alpha}  \right)\\
&&\qquad \qquad \qquad \qquad \qquad \qquad \times 
\left(\zeta(1-\gamma+\delta) u^{-\gamma}q^{\gamma-\delta}
+  \zeta(1+\gamma-\delta) u^{-\delta}q^{\delta-\gamma}  \right).
\end{eqnarray*}
We multiply out this expression and have 4 terms each of  which can be expressed in terms of  
\begin{eqnarray*}
M_{\alpha,\beta;\gamma,\delta}:  
=\frac{\zeta(1+\alpha+\beta)\zeta(1+\gamma+\delta)}{\zeta(2+\alpha+\beta+\gamma)}
\sum_{d\mid h} \frac{1}{d^{1+\alpha+\beta+\gamma+\delta}}.
\end{eqnarray*} 
Specifically, we have 
$$M'(u,h)= u^{-\alpha-\gamma} M_{-\gamma,\beta;-\alpha,\delta}+ u^{-\alpha-\delta}M_{-\delta,\beta;\gamma,-\alpha}+
 u^{-\beta-\gamma}M_{\alpha,-\gamma; -\beta,\delta}+ u^{-\beta-\delta}M_{\alpha,-\delta;\gamma,-\beta}
.$$
Now, much as we did in the first section we compute 
$$2\Re\int_0^\infty 
\sum_{h\le v/y} M'(\frac{hT}{2\pi v},h) \hat{\psi}(v) \frac{dv}{v}.$$
(Note that the integrand is 0 for $0<v<y$.)
We have that one term (the term corresponding to $u^{-\alpha-\gamma} M_{-\gamma,\beta;-\alpha,\delta}$)
of $\sum_{h\le v/y} M'(\frac{hT}{2\pi v},h) $ is  
\begin{eqnarray*}&&
\frac{\zeta(1-\alpha+\beta)\zeta(1-\gamma+\delta)}{\zeta(2-\alpha+\beta-\gamma+\delta)}\left(\frac{T}{2\pi v}\right)^{-\alpha-\gamma}
\sum_{hd\le v/y} h^{-\alpha-\gamma}d^{-1 -\beta -\delta}\\
&&  = \frac{\zeta(1-\alpha+\beta)\zeta(1-\gamma+\delta)}{\zeta(2-\alpha+\beta-\gamma+\delta)}\left(\frac{T}{2\pi v}\right)^{-\alpha-\gamma}
\frac{1}{2\pi i}\int_{(2)}  \zeta(s+\alpha+\gamma) \zeta(s+1 \beta +\delta)\frac{ (\frac vy  )^s}{s}~ds.
\end{eqnarray*}
As before, see (\ref{eqn:mellin}), we compute that  
\begin{eqnarray*}
\zeta(s+\alpha+\gamma) \int_0^\infty \hat{\psi}(v) v^{s-1+\alpha+\gamma} ~dv = \zeta(1-s-\alpha-\gamma)
\int_0^\infty \psi(u) u^{-s-\alpha-\gamma} ~du.
\end{eqnarray*}
Then, by moving the $s$-path of integration to $\Re s=-1/2$  we have 
\begin{eqnarray*}&&
\frac{1}{2\pi i} \int_{(2)} \zeta(1-s-\alpha-\gamma)\zeta(s+1 +\beta +\delta)\frac{y^{-s}u^{-s-\alpha-\gamma}}{s}~ds\\
&&
= u^{-\alpha-\gamma}\zeta(1-\alpha-\gamma)\zeta(1 +\beta +\delta)+
 \zeta(1-\alpha +\beta-\gamma +\delta)\frac{ y  ^{\alpha+\gamma}}{\alpha+\gamma}\\
&& \qquad \qquad \qquad 
-\zeta(1-\alpha+\beta-\gamma+\delta)  \frac{  y^{\beta+\delta}u^{-\alpha+ \beta -\gamma+\delta}}{ \beta +\delta} +O(y^{1/2}).
\end{eqnarray*}
Now we have shown that 
\begin{eqnarray*}&&
2\Re\int_0^\infty 
\sum_{h\le v/y} M'(\frac{hT}{2\psi v},h) \hat{\psi}(v) \frac{dv}{v}=\int_0^\infty \psi(u)
\left(\frac{T}{2\pi }\right)^{-\alpha-\gamma}\frac{\zeta(1-\alpha+\beta)\zeta(1-\gamma+\delta)}{\zeta(2-\alpha+\beta-\gamma+\delta)}\\
&&\qquad \times 
\bigg(u^{-\alpha-\gamma}\zeta(1-\alpha-\gamma)\zeta(1 +\beta +\delta)+
 \zeta(1-\alpha +\beta-\gamma +\delta)\frac{ y  ^{\alpha+\gamma}}{\alpha+\gamma}\\
&& \qquad \qquad \qquad 
-\zeta(1-\alpha+\beta-\gamma+\delta)  \frac{  y^{\beta+\delta}u^{-\alpha+ \beta -\gamma+\delta}}{ \beta +\delta}\bigg) ~du +O(y^{1/2}).
\end{eqnarray*}
Recalling that $y=T/(2\pi X)$ we see that
this term matches up with the first term (\ref{eqn:oneswap}) that we got from the recipe with one swap.
Similarly the other three terms match up.  
It follows that the one-swap terms match up with the Type I off-diagonal terms. 

\section{Conclusion}

The coefficient correlation terms that we have analyzed here contribute to $I(T;X)$ once $X$ gets larger than $T$.
We have shown two very different ways to view these terms: from the 1-swap terms of the recipe conjecture,
and from the Type-I divisor correlation sums much studied in analytic number theory. 
Our results illustrate that there are other terms to be considered  once $X>T^2$, i.e. that there is something missing from the above analysis.  These are the divisor correlations that must match the 2-swapped terms.

The next step is to describe, from the ascending perspective,  the terms  
that  we see in the  descent   from (\ref{eqn:fourth}) to  the terms with 2-swaps. This will be the subject of a subsequent paper where
the coefficient sums to be considered are much as in [BK]:
\begin{eqnarray*}\sum_{M_1,M_2,h_1,h_2\atop(M_1,M_2)=1}
\sum_{m_1 , m_2,n_1,n_2 \atop {(*_1), (*_2)\atop m_1m_2\le X}}\frac{ m_1^{-\alpha}   n_1^{-\gamma} m_2^{-\beta}
 n_2^{-\delta}}{m_1m_2}
\hat \psi\left( T\left(\frac{h_1}{m_1M_1}+\frac{h_2}{m_2M_2}\right) \right)
\end{eqnarray*}
where 
$$(*_1): M_1 m_1=M_2n_1+h_1 \qquad \qquad (*_2): M_2 m_2=M_1n_2+h_2$$
for coprime integers $M_1$ and $M_2$.

 \end{document}